\documentclass[letterpaper, 10 pt, conference]{ieeeconf}
\usepackage{blkarray}
\usepackage{dsfont}
\usepackage[cmex10]{amsmath}
\usepackage{bm}
\usepackage{comment}
\usepackage{amssymb}
\usepackage{gensymb}
\usepackage{epsfig,floatflt}
\usepackage{psfrag,graphicx}
\usepackage{algorithm}
\usepackage{algpseudocode}
\usepackage[colorlinks]{hyperref}
\usepackage{cite}
\usepackage{epstopdf}
\usepackage{subcaption}
\usepackage[labelformat=default,labelsep=colon]{caption}
\IEEEoverridecommandlockouts

\newcommand{\Input}{\item[\textbf{Input:}]}
\newcommand{\Output}{\item[\textbf{Output:}]}

\newcommand{\IR}{\mathbb{R}}

\newtheorem{prob}{\textbf{Problem}}

\newtheorem{rem}{\textbf{Remark}}

\DeclareSymbolFont{extraup}{U}{zavm}{m}{n}
\DeclareMathSymbol{\vardiamond}{\mathalpha}{extraup}{87}
\newcommand*{\qedprob}{\hfill\ensuremath{\triangledown}}%
\newcommand*{\qedproba}{\ensuremath{\triangledown}}%

\begin{document}
	

    \title{Real-Time Optimal Control via Transformer Networks and Bernstein Polynomials}

	\author{Gage MacLin$^{1}$, Venanzio Cichella$^{1}$, Andrew Patterson$^{2}$, and Irene Gregory$^{2}$
		\thanks{This work was supported by NASA and ONR.}
		\thanks{$^{1}$ Gage MacLin and Venanzio Cichella are with the Department of Mechanical Engineering, University of Iowa, Iowa City, IA 52240 {\tt\small \{gage-maclin, venanzio-cichella\} @uiowa.edu}}
        \thanks{$^{2}$ Andrew Patterson and Irene Gregory are with the NASA Langley Research Center, Hampton, Virginia, VA 23681 {\tt\small \{andrew.patterson, irene.m.gregory\} @nasa.gov}}
	}

	\maketitle

\begin{abstract}
 In this paper, we propose a Transformer-based framework for approximating solutions to infinite-dimensional optimization problems: calculus of variations problems and optimal control problems. Our approach leverages offline training on data generated by solving a sample of infinite-dimensional optimization problems using composite Bernstein collocation. Once trained, the Transformer efficiently generates near-optimal, feasible trajectories, making it well-suited for real-time applications. In motion planning for autonomous vehicles, for instance, these trajectories can serve to warm-start optimal motion planners or undergo rigorous evaluation to ensure safety. We demonstrate the effectiveness of this method through numerical results on a classical control problem and an online obstacle avoidance task. This data-driven approach offers a promising solution for real-time optimal control of nonlinear, nonconvex systems.
\end{abstract}


	\IEEEpeerreviewmaketitle


\section{Introduction}\label{sec:Introduction}
Solving optimal control problems (OCPs) in real time poses significant computational challenges, especially in safety-critical applications such as robotics, autonomous vehicle navigation, and spacecraft guidance. These problems often involve nonlinear dynamics and nonconvex constraints, making them computationally expensive to solve. In real-world scenarios, the ability to generate feasible trajectories in real time is crucial for collision avoidance, dynamic replanning, and multi-agent coordination. Model predictive control (MPC) and convex optimization are two widely used frameworks for online trajectory generation. MPC offers a powerful framework for dynamically updating control inputs based on real-time feedback while handling constraints. Convex optimization methods, on the other hand, provide efficient and reliable solutions by exploiting problem convexity. Despite significant progress in extending these methods \cite{grune2016nonlinear,malyuta2022convex}, they struggle to maintain real-time capabilities for more complex problems. This has motivated the exploration of data-driven methods for the online solution of OCPs.

Transformers \cite{vaswani2017attention}, originally developed for natural language processing tasks, have emerged as a powerful architecture for sequence modeling due to their ability to learn long-range dependencies in sequential data. Unlike traditional sequence modeling methods, such as recurrent neural networks (RNNs) and long short-term memory (LSTM) networks, which both process data sequentially, Transformers handle entire sequences in parallel. This parallelism can significantly boost efficiency during training, especially when processing long sequences. Additionally, their flexibility and capacity to capture complex input-output relationships makes them highly effective for modeling intricate sequential data. These properties make Transformers well suited for problems dependent on understanding complex temporal relationships, such as trajectory generation \cite{maclin2025transformer}. 

The potential of Transformers for solving OCPs lies in their ability to learn the underlying structure of optimal trajectories from a dataset of precomputed solutions. When trained offline on a sufficiently large and diverse set of OCP solutions, a Transformer can generalize to new scenarios and generate potentially feasible trajectories online, avoiding the need to solve the OCP from scratch. Compared to traditional numerical methods, this data-driven approach offers substantially greater computational efficiency. Alternative machine learning-based methods, e.g., \cite{hauser2017learning}, also leverage precomputed solutions to generate trajectories online, though they are typically limited in generalization and scalability. In contrast, Transformers excel at capturing complex, long-range dependencies across the solution space, enabling them to generate potentially feasible trajectories even for scenarios not explicitly seen during training.

However, a major challenge in using data-driven methods for real-time optimal control is that they lack guarantees on constraint satisfaction. Neural networks, including Transformers, offer no inherent guarantee that the generated solutions satisfy the system dynamics or problem constraints. To mitigate this, the predicted trajectories can be used to warm-start an OCP \cite{guffanti2024transformers,celestini2024transformer}. By initializing the solver with a high-quality guess, the number of solver iterations and overall computation time can be reduced substantially compared to a cold start, leading to noticeably faster convergence in practice. Alternatively, if the network output is designed to follow a structured form that inherently satisfies constraints and safety conditions, this structure can be leveraged for efficient verification, e.g., \cite{patterson2023hermite}. In this case, feasibility and safety can be confirmed without the need for a post-inference solver, enhancing the reliability of the data-driven solution.


In this work, we propose a Transformer-based framework for real-time optimal control. The network takes as input a set of parameters that define a problem's constraints and objective function and outputs the solution parameterized as coefficients of composite Bernstein polynomials. This parameterization leverages the geometric properties of Bernstein polynomials \cite{kielas2022bernstein,cichella2020optimal}, enabling rigorous validation of the neural network's solutions. To train the model, we first generate a large dataset of optimal solutions by solving problems with sampled parameter values using composite Bernstein collocation \cite{maclin2024optimal}. The Transformer then learns the underlying patterns of these solutions. Once trained, it can generate potentially feasible trajectories online. Moreover, the structured output, expressed as coefficients of composite Bernstein polynomials, facilitates feasibility verification. Specifically, the degree elevation, convex hull, and endpoint properties of Bernstein polynomials \cite{kielas2022bernstein} can be exploited to ensure constraint satisfaction. This hybrid approach significantly reduces the computational burden of solving OCPs in real time while preserving feasibility guarantees.



This paper is organized as follows: in Section \ref{sec:OCP} we present two classes of infinite-dimensional optimization problems. Section \ref{sec:CBP} introduces composite Bernstein polynomials and highlights some of their properties. Section \ref{sec:NLP} presents the discretized OCP. In Section \ref{sec:TransformersOCP}, the Transformer framework for solving OCPs is introduced. Numerical examples are discussed in Section \ref{sec:Num}, including the Brachistochrone problem and an online motion planning problem. Then, conclusions are presented in Section \ref{sec:Conclusions}.

\section{Infinite-Dimensional Optimization Problems}\label{sec:OCP}
In this work, we focus on two classes of infinite-dimensional optimization problems. First, we consider the problem of calculus of variations:

\begin{prob}[Problem $P^{CV}$]\label{prob:calculusofvariations}

Determine $\bm{x}(t):[0,t_f]\to\IR^{n_x}$ and $t_f$ that minimizes the modified Bolza-type cost functional
	\begin{equation} \label{eq:costfunc2}
	\min_{\bm{x}(t),{{t}_f}}  J(\bm{x}(t),{{t}_f})=
    E(\bm{x}(0),\bm{x}(t_f),t_f)+\int^{t_f}_0F(\bm{x}(t))dt, 
	\end{equation}
subject to equality constraints
\begin{equation}\label{eq:equalityconstraint2}
    \bm{e}(\bm{x}(0),\bm{x}(t_f),{{t}_f}) = \bm{0}, 
\end{equation}
and inequality constraints
\begin{equation}\label{eq:inequalityconstraint2}
    \bm{h}(\bm{x}(t)) \leq \bm{0} \, , \quad \forall t\in [0,t_f] ,
\end{equation}
where $J:\mathbb{R}^{n_x}\times\mathbb{R}\to\mathbb{R}$, $E:\IR^{n_x}\times\IR^{n_x}\times\IR\to\IR$, $F:\IR^{n_x}\to\IR$, $\bm{e}:\IR^{n_x}\times\IR^{n_x}\times\IR\to\IR^{n_e}$, and $\bm{h}:\IR^{n_x}\to\IR^{n_h}$.\qedproba
\end{prob}

Then, with the introduction of system dynamics, the following optimal control problem is formulated:
\begin{prob}[Problem $P^{OC}$]\label{prob:continuous}

Determine $\bm{x}(t):[0,t_f]\to\IR^{n_x}$, $\bm{u}(t):[0,t_f]\to\IR^{n_u}$ and $t_f$ that minimizes the Bolza-type cost functional
	\begin{equation} \label{eq:costfunc}
	\begin{split}  
	& \min_{\bm{x}(t),\bm{u}(t),{{t}_f}}  \tilde{J}(\bm{x}(t),\bm{u}(t),{{t}_f}) = \\  &  \tilde{E}(\bm{x}(0),\bm{x} 
 (t_f),t_f)+ \int_0^{t_f} \tilde{F}(\bm{x}(t),\bm{u}(t))dt \, ,
	\end{split} 
	\end{equation}
subject to dynamic constraints
\begin{equation}\label{eq:dynamicconstraint}
    \dot{\bm{x}}(t) = \bm{f}(\bm{x}(t),\bm{u}(t)) , \quad \forall t \in[0,t_f],
\end{equation}
equality constraints
\begin{equation}\label{eq:equalityconstraint}
    \tilde{\bm{e}}(\bm{x}(0),\bm{x}(t_f),{{t}_f}) = \bm{0}, 
\end{equation}
and inequality constraints
\begin{equation}\label{eq:inequalityconstraint}
    \tilde{\bm{h}}(\bm{x}(t),\bm{u}(t)) \leq \bm{0} \, , \quad \forall t\in [0,t_f] ,
\end{equation}
where $\tilde{J}:\IR^{n_x} \times \IR^{n_u} \times \IR \to \IR$, $\tilde{E}:\IR^{n_x}\times \IR^{n_x} \times \IR \to \IR$, $\tilde{F}:\IR^{n_x} \times \IR^{n_u} \to \IR$, $\bm{f}:\IR^{n_x}\times \IR^{n_u}\to \IR^{n_x}$, $\tilde{\bm{e}}:\IR^{n_x} \times \IR^{n_x} \times \IR \to \IR^{n_e}$, and $\tilde{\bm{h}}:\IR^{n_x} \times \IR^{n_u}  \to \IR^{n_h}$.
	\qedprob
\end{prob}


Then, we define a problem space by considering variations in the functions \(E\), \(F\), \(\bm{e}\), and \(\bm{h}\), encapsulated through a problem-specific parameter denoted as \(\bm{\theta} \in \Theta \). To explicitly express the dependence of the problem on \(\bm{\theta}\), and with a slight abuse of notation, we rewrite Problem $P^{CV}$ as follows: 

\begin{prob}[Problem $P^{CV}(\bm{\theta})$]\label{prob:calculusofvariations_theta}

Determine $\bm{x}(t)$ and $t_f$ that minimizes
	\begin{equation} \label{eq:costfunc4}
    \begin{split}
        & \min_{\bm{x}(t),{{t}_f}}  J(\bm{x}(t),{{t}_f},\bm{\theta}) =   \\
        & E(\bm{x}(0),\bm{x}
        (t_f),\bm{\theta})+ \int_0^{t_f} F(\bm{x}(t),\bm{\theta})dt \, ,
    \end{split}
	\end{equation}
\begin{equation}\label{eq:equalityconstraint4}
    \bm{e}(\bm{x}(0),\bm{x}(t_f),{{t}_f},\bm{\theta}) = \bm{0}, 
\end{equation}
\begin{equation}\label{eq:inequalityconstraint4}
    \bm{h}(\bm{x}(t),\bm{\theta}) \leq \bm{0} \, , \quad \forall t\in [0,t_f] ,
\end{equation}
where $E(\bm{x}(0),\bm{x}(t_f),\bm{\theta})=E(\bm{x}(0),\bm{x}(t_f))$, $F(\bm{x}(t),\bm{\theta})=F(\bm{x}(t))$, $\bm{e}(\bm{x}(0),\bm{x}(t_f),t_f,\bm{\theta})=\bm{e}(\bm{x}(0),\bm{x}(t_f),t_f)$, and $\bm{h}(\bm{x}(t),\bm{\theta})=\bm{h}(\bm{x}(t))$.
    \qedprob
\end{prob}

The above reparameterization can be similarly applied to Problem $P^{OC}$ to obtain Problem $P^{OC}(\bm{\theta})$.

\section{Composite Bernstein Approximation}\label{sec:CBP}
In this section, the composite Bernstein polynomial will be presented, along with some discussion of their beneficial properties for the solution of optimal control and calculus of variations problems. A composite Bernstein polynomial is a polynomial that consists of $K$ Bernstein polynomials of degree $N$ connected in series at time knots $t_k$. An $N$th order Bernstein polynomial is defined over the domain $I_k=t\in[t_{k-1},t_k]$ as
\begin{equation}
    x^{[k]}_N(t)=\sum^N_{j=0}\bar{x}_{j,N}^{[k]}b_{j,N}^{[k]}(t),
\end{equation}
where $k$ represents the subdomain of the Bernstein polynomial, $\bar{x}_{j,N}^{[k]}$ are the Bernstein control points, and $b_{j,N}^{[k]}$ is the Bernstein basis function
\begin{equation}
    b_{j,N}^{[k]}(t)=\binom{N}{j}\frac{(t-t_{k-1})^j(t_k-t)^{N-j}}{(t_k-t_{k-1})^N},
\end{equation}
for $j=0,...,N$. Now, we can define a composite Bernstein polynomial as
\begin{equation}
    x_M(t)=x_N^{[k]}(t), \quad \forall t \in I_k, \quad\forall k=1,...,K.
\end{equation}
With each Bernstein polynomial defined over a subdomain of $I_k$, we define the time knots $t_k$ of a composite Bernstein polynomial to be the extremes of each subdomain $[t_{k-1},t_k]$, where $k=0,...,K$, with the assertion that $t_0<t_1<...<t_K$. The control points are then defined as
\begin{equation}\label{eq:CBPCP}
    \bar{x}_M=\left[\bar{x}_{0,N}^{[1]},...,\bar{x}_{N,N}^{[1]},...,\bar{x}_{0,N}^{[K]},...,\bar{x}_{N,N}^{[K]}\right],
\end{equation}
where $\bar{x}_M$ is said to have $M+1$ control points with $M=K(N+1)-1$, ensuring consistency with the notation of standard Bernstein polynomials.

The derivative of a composite Bernstein polynomial is defined over the interval $[t_0,t_K]$ as
\begin{equation}\label{eq:CBPderivative}
    \dot{x}_M(t)=\sum^{N-1}_{j=0}\sum^N_{i=0}\bar{x}_{i,N}^{[k]}D_{i,j}^{[k]}b_{j,N-1}^{[k]}(t),
\end{equation}
where $D_{i,j}^{[k]}$ is the $(i,j)$th entry of the Bernstein differentiation matrix
\begin{multline}
    \bm{D}_{N-1}^{[k]}=
    \begin{bmatrix}
- \frac{N}{t_k-t_{k-1}}  &  0 & \ldots & 0 \\
\frac{N}{t_k-t_{k-1}} &  \ddots & \ddots & \vdots \\ 
0 &  \ddots & \ddots & 0 \\
\vdots & \ddots & \ddots  & -\frac{N}{t_k-t_{k-1}} \\
0 &  \ldots & 0 & \frac{N}{t_k-t_{k-1}}
\end{bmatrix}\\
\in\mathbb{R}^{(N+1)\times N}, \; \forall k=1,...,K.
\end{multline}

The definite integral of a composite Bernstein polynomial is obtained via
\begin{equation}
    \int^{t_K}_{t_0}x_M(t)dt=\sum^K_{k=1}w^{[k]}\sum^N_{j=0}\bar{x}_{j,N}^{[k]},
\end{equation}
where 
\begin{equation}\label{eq:CBPweights}
    w^{[k]}=\frac{t_k-t_{k-1}}{N+1} , \; \forall k = 1, \ldots, K.    
\end{equation}

Beyond being differentiable and integrable, composite Bernstein polynomials have two particularly beneficial properties for the solution of optimal control problems. First being the end point property, which stipulates that the terminal values of each polynomial are equivalent to the terminal control points, e.g.,
\begin{equation}
    x_M(t_k)=\bar{x}_{N,N}^{[k-1]}=\bar{x}_{0,N}^{[k]},
\end{equation}
which implicitly enforces boundary constraints. Secondly, composite Bernstein polynomials have the convex hull property, which states that the polynomial lies entirely within the convex hull defined by its control points. This property offers a mathematical tool for guaranteeing satisfaction of inequality constraints for the entire polynomial, not just at the control points. Notably, the resulting guarantee can be quite conservative \cite{tordesillas2022minvo}; however, we can refine the convex hull to get a tighter representation via degree elevation, i.e.,

\begin{equation}
    \bar{x}^{[k]}_{j,N+1} = \bar{x}^{[k]}_{j,N}\bm{E}^{N_e}_N, \quad \forall k=1,...,K,
\end{equation}
where $\bm{E}^{N_e}_N=\{e_{i,j}\}\in\IR^{(N+1)\times(N_e+1)}$ is the degree elevation matrix from degree $N$ to $N_e>N$, with its elements defined by
\begin{equation}
    e_{i,i+j}=\begin{cases}
        \frac{\binom{N_e-N}{j}\binom{N}{i}}{\binom{N_e}{i+j}}, & i\le N, \,\,\,j\le N_e-N, \\
        \quad\,\,\,\,0, & \text{otherwise.}
    \end{cases}
\end{equation}
Moreover, composite Bernstein polynomials yield consistent approximations of continuous optimal control problems, they also achieve a higher rate of convergence to the solution compared to a single Bernstein approximant \cite{maclin2024optimal, cichella2019consistent}. To this end, we use composite Bernstein collocation to discretize Problems $P^{CV}(\bm\theta)$ and $P^{OC}(\bm\theta)$.

\section{Direct Approximation of Infinite-Dimensional Optimization Problems}\label{sec:NLP}

The output of Problem $P^{CV}(\bm\theta)$ can be approximated using a composite Bernstein polynomial $\bm{x}_M(t):[t_0,t_K]\to\IR^{n_x}$. That is, $\lim_{M\to\infty}\bm{x}_M(t)=\bm{x}(t)$,
where
\begin{equation}\label{eq:CBPx}
    \bm{x}_M(t)=\sum^N_{j=0}\bar{\bm{x}}_{j,N}^{[k]}b_{j,N}(t), \quad t\in[t_{k-1},t_k], \,\, \forall k=1,...,K.
\end{equation}
Letting $\bar{\bm{x}}_M$ be the control points of $\bm{x}_M(t)$ and defined in the same form as Equation \eqref{eq:CBPCP}, we can formally define the discretized Problem $P^{CV}(\bm{\theta})$ as:
\begin{prob}[Problem $P_M^{CV}(\bm{\theta})$]\label{prob:discrete}
Determine $\bar{\bm x}_M$ and $t_K$ that minimizes
\begin{equation} \label{eq:costfunc3}
	\begin{split}  
	& \min_{\bar{\bm{x}}_M,{{t}_K}}  J(\bar{\bm{x}}_M,{{t}_K},\bm{\theta}) = \\  &  E(\bar{\bm{x}}_{0,M},\bar{\bm{x}}_{M,M},t_K,\bm{\theta})+ \sum^K_{k=1}w^{[k]}\sum^N_{j=0}F(\bar{\bm{x}}_{j,N}^{[k]},\bm{\theta}),
	\end{split} 
	\end{equation}
subject to
\begin{equation}\label{eq:equalityconstraint3}
    \bm{e}(\bar{\bm{x}}_{0,M},\bar{\bm{x}}_{M,M},{{t}_K},\bm{\theta}) = \bm{0}, 
\end{equation}
\begin{equation}\label{eq:inequalityconstraint3}
    \bm{h}(\bar{\bm{x}}_{j,M},\bm{\theta}) \leq \bm{0} \, , \quad \forall j=0,...,M,
\end{equation}
\begin{equation}\label{eq:knotcontinuity}
\bar{\bm{x}}_{N,N}^{[k]}-\bar{\bm{x}}_{0,N}^{[k+1]}=0, \forall k=1,...,K-1,
\end{equation}
\begin{equation}\label{eq:timecontinuity}
    t_K>t_{K-1}>\dots>t_0,
\end{equation}
where $w^{[k]}$ is defined in Equation \eqref{eq:CBPweights}
\qedprob
\end{prob}

The output of Problem $P^{OC}(\bm\theta)$ can be similarly approximated using composite Bernstein polynomials $\bm x_M(t):[t_0,t_k]\to\IR^{n_x}$ and $\bm u_M(t):[t_0,t_K]\to\IR^{n_u}$, defined in the same form as Equation \eqref{eq:CBPx}. Then, with $\bar{\bm u}_M$ defined as the control points of $\bm u_M(t)$, the discretized Problem $P^{OC}(\bm\theta)$ follows:

\begin{prob}[Problem $P_M^{OC}(\bm{\theta})$]\label{prob:discrete_oc}
Determine $\bar{\bm x}_M$, $\bar{\bm u}_M$ and $t_K$ that minimizes
\begin{equation} \label{eq:costfunc3}
	\begin{split}  
	& \min_{\bar{\bm{x}}_M,\bar{\bm u}_M,{{t}_K}}  J(\bar{\bm{x}}_M,\bar{\bm u}_M,{{t}_K},\bm{\theta}) = \\  &  E(\bar{\bm{x}}_{0,M},\bar{\bm{x}}_{M,M},t_K,\bm{\theta})+ \sum^K_{k=1}w^{[k]}\sum^N_{j=0}F(\bar{\bm{x}}_{j,N}^{[k]},\bar{\bm{u}}_{j,N}^{[k]},\bm{\theta}),
	\end{split} 
	\end{equation}
subject to
\begin{equation}\label{eq:dynamicconstraint3}
\left\Vert\sum^M_{i=0}\bar{\bm x}_{i,M}\bm{D}_{i,j}^M-\bm{f}(\bar{\bm x}_{j,M},\bar{\bm u}_{j,M},\bm\theta)\right\Vert\le\delta_P^M, \,\forall j=0,...,M
\end{equation}
\begin{equation}\label{eq:equalityconstraint3}
    \bm{e}(\bar{\bm{x}}_{0,M},\bar{\bm{x}}_{M,M},{{t}_K},\bm{\theta}) = \bm{0}, 
\end{equation}
\begin{equation}\label{eq:inequalityconstraint3}
    \bm{h}(\bar{\bm{x}}_{j,M},\bar{\bm{u}}_{j,M},\bm{\theta}) \leq \bm{0} \, , \quad \forall j=0,...,M,
\end{equation}
\begin{equation}\label{eq:knotcontinuity}
\bar{\bm{x}}_{N,N}^{[k]}-\bar{\bm{x}}_{0,N}^{[k+1]}=0, \forall k=1,...,K-1,
\end{equation}
\begin{equation}\label{eq:timecontinuity}
    t_K>t_{K-1}>\dots>t_0,
\end{equation}
where $\delta_P^M$ is a relaxation bound that converges uniformly to 0, and $\bm{D}_{i,j}^M$ is the $(i,j)$th entry of the differentiation matrix described in \cite{maclin2024optimal}.
\qedprob
\end{prob}

\begin{rem}
    The solution to Problem $P_M^{CV}(\bm{\theta})$ provides a set of optimal control points $\bar{\bm{x}}^*_M$ that yields the optimal trajectory $\bm{x}_M^*(t)$. Proofs establishing the feasibility and convergence of Problem $P^{CV}_N$ to Problem $P^{CV}$, i.e., $\lim_{N\to\infty}\bm{x}_N^*(t)=\bm{x}^*(t)$ are provided in \cite{cichella2017optimal} and extended to composite Bernstein collocation in \cite{maclin2024optimal}. This remark holds for the solution to Problem $P^{OC}_M(\bm\theta)$.
\end{rem}

\section{Transformers}\label{sec:TransformersOCP}
A Transformer is a neural network consisting of an encoder and a decoder, with additional pre- and post-processing. The encoder processes the input data sequence and creates a high dimensional representation that can capture the relationships between each data point. Then, based on this high dimensional representation, the decoder generates the output data sequence one value at a time. 

\begin{figure}[H]
    \centering
    \includegraphics[width=.48\textwidth]{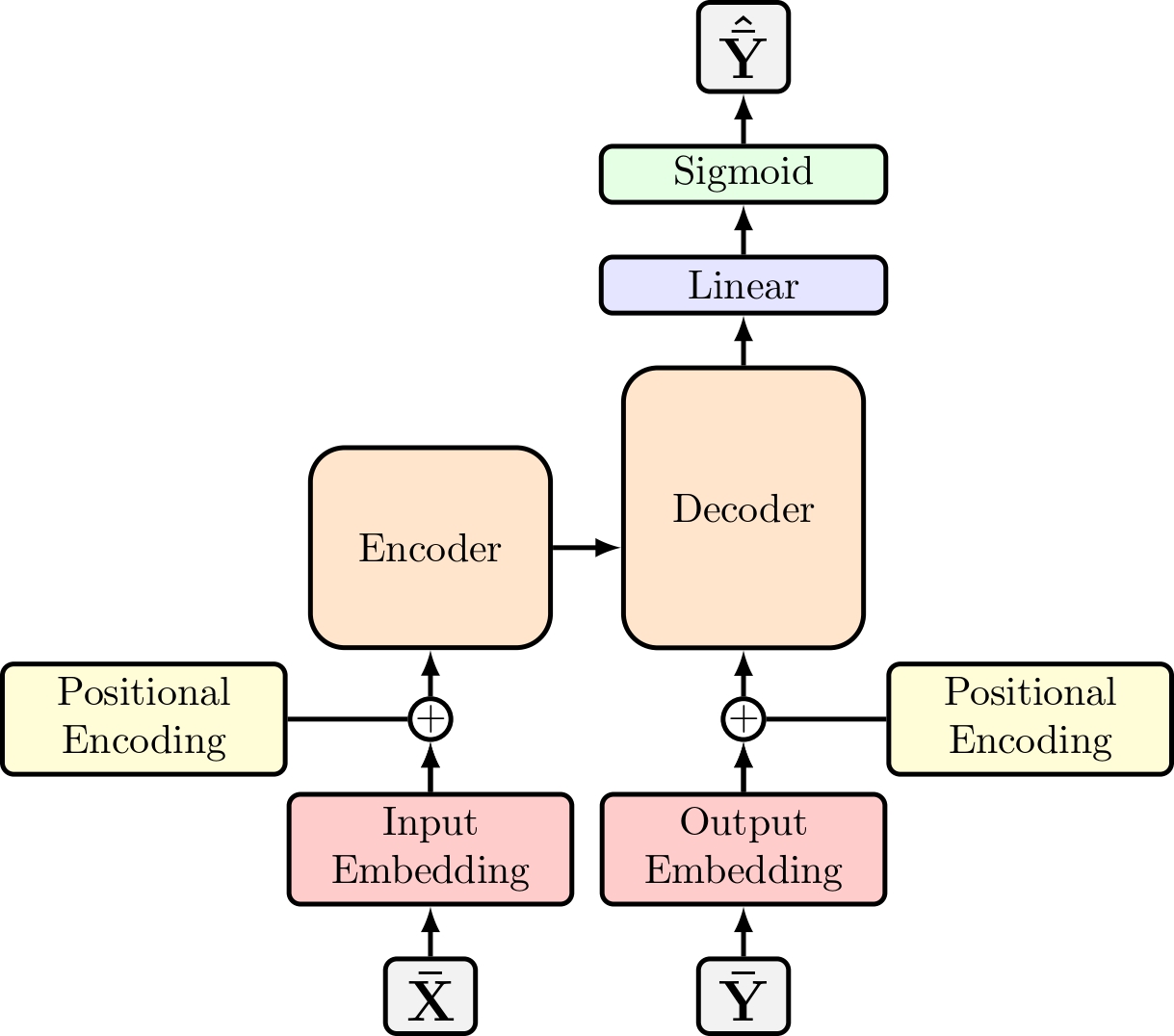}
    \caption{Transformer architecture: Input state $\bm{\bar{X}}$ and output state $\bm{\bar{Y}}$ are pre-processed and provided to the encoder and decoder respectively, which yields the predicted output state $\bm{\hat{\bar{Y}}}$ after post-processing.}
    \label{fig:transformer}
\end{figure}

In this work, we use the standard encoder-decoder framework from PyTorch \cite{paszke2019pytorch}. We explicitly tailor the aforementioned pre- and post-processing to create a Transformer model that can effectively predict solutions to Problem $P_M^{CV}(\bm{\theta})$. This architecture is shown in Figure \ref{fig:transformer}. First the input state $\bm{\bar{X}}$ and output state $\bm{\bar{Y}}$  are each embedded. This translates the data into a higher dimension, allowing for the Transformer to build a higher level understanding of the data. The embedded $\bm{\bar{X}}$ and $\bm{\bar{Y}}$ are then combined with learned positional encodings, which adaptively capture the sequential nature of the data, before being passed to the encoder and decoder, respectively. The decoder output is then processed back to an appropriate dimension (linear layer), and then scaled s.t. $\bm{\hat{\bar{Y}}}\in[0,1]$ (sigmoid layer).



\subsection{Training}
With the objective of generating solutions to optimal control problems online, we must first generate $N_d$ data sequences offline by solving a discretized infinite-dimensional optimization problem for various $\bm{\theta}$. That is, we solve $P^{CV}_M(\bm{\theta}^i)$ or $P^{OC}_M(\bm{\theta}^i)$ for $i=1,...,N_d$ where $\bm{\theta}^i$ are randomly sampled, to create a diverse dataset. Once we have this dataset, we can train a Transformer model to solve the same problem online. To this end, we design a teacher forcing training algorithm that uses a causal self-attention mask to learn the dependencies between sequential data points. Then, we minimize the squared error loss function:
\begin{equation}\label{eq:loss}
    \mathcal{L}=\sum_{n=1}^{N_d}\sum_{i=0}^M||\bar{\bm{y}}_{i,M}^n-\hat{\bar{\bm{y}}}_{i,M}^n||_2^2,
\end{equation}
where $\hat{\bar{\bm{y}}}_{i,M}^{n}$ are the Transformer-predicted control points, and $\bar{\bm{y}}_{i,M}^{n}$ are the known optimal control points for each sequence $n$ in the dataset. For problems of type $P^{CV}_M(\bm\theta)$, $\bar{\bm y}_M=\bar{\bm x}_M$, whereas for problems of type $P^{OC}_M(\bm\theta)$, $\bar{\bm y}_M=[\bar{\bm x}_M,\bar{\bm u}_M]$.

\subsection{Inference}
We design a general inference pipeline to solve infinite-dimensional optimization problems using composite Bernstein polynomials. Specifically, we use a trained Transformer model to predict the control points $\bar{\bm y}_M$. The model takes as input the problem-specific parameter $\bm{\theta}$, and generates a complete set of control points through autoregression (see Algorithm \ref{alg:inference}). This framework is inherently flexible and can be adapted to problem-specific requirements. 


\begin{algorithm}[t]
\caption{Inference}\label{alg:inference}
\begin{algorithmic}[1]
\Input $\bm{\theta}$
\Output $\bar{\bm{y}}_{i,M}, \forall i\in\{0,M\}$
\State $\bar{\bm{X}}=\bm{\theta}$
\For{$i = 0$ to $M-1$}
    \State $\bar{\bm{Y}}=\bar{\bm{y}}_{j,M}, \forall j\in\{0,i\}$
    \State $\bar{\bm{y}}_{i+1,M}=\text{Transformer}(\bar{\bm{X}},\bar{\bm{Y}})$
\EndFor
\State
\Return $\bar{\bm{y}}_{i,M}, \forall i\in\{0,M\}$
\end{algorithmic}
\end{algorithm}


\subsection{Feasibility of trajectories}
Using a Transformer provides no guarantees on constraint satisfaction. However, with the data consisting of the control points of composite Bernstein polynomials instead of just sampled data points, feasibility can be quickly verified due to the properties discussed in Section \ref{sec:CBP}, namely the degree elevation, convex hull, and end point properties \cite{kielas2022bernstein}. Alternatively, the outputs of the Transformer model can be used to warm-start the original optimal control problem, potentially improving the efficiency of the nonlinear program significantly.

\section{Numerical Results}\label{sec:Num}
\subsection{Brachistochrone problem}
To investigate the efficacy of this approach, we turn to the classic Brachistochrone problem:
\begin{equation}
    \min_{y(x)} J(y(x))=\int^{\theta_1}_0\frac{\sqrt{1+\dot{y}(x)^2}}{\sqrt{2gy(x)}}dx,
\end{equation}
subject to equality constraints
\begin{equation}
    y(0)=0,\quad y(\theta_1)=\theta_2,
\end{equation}
where g is the gravitational constant, and $\theta_1$ and $\theta_2$ are the final $x$ and $y$ position respectively. Conveniently, the Brachistochrone problem is a calculus of variations problem with a known analytical solution \cite{liberzon2011calculus}. We construct a sufficiently large dataset of $N_d=100,000$ samples by computing the analytical solution of the problem with different values of $\bm{\theta}^i$, where $i=1,...,N_d$, and parameterize it using composite Bernstein polynomials. The coefficients of these polynomials serve as the dataset used to train the Transformer, using the hyperparameters shown in Table \ref{tab:parameters}.

\begin{table}
    \centering
    \caption{Hyperparameter configuration.}
    \label{tab:parameters}
    \begin{tabular}{|l|c|c|}
        \hline
        \textbf{Hyperparameter}  & \textbf{Transformer} \\ 
        \hline
        Learning rate                   & $4 \times 10^{-4}$ \\
        Batch size                    & $8$ \\
        Number of layers                 & $1$  \\
        Embedding dimension $(d_{\text{emb}})$               & $64$  \\
        Dropout rate                    & $0.1$ \\
        Number of attention heads          & $4$ \\
        Optimizer                            & Adam  \\
        Dataset size $(N_d)$        & $100,000$ \\
        \hline
    \end{tabular}
\end{table}

\begin{figure}
    \centering
    \includegraphics[width=.48\textwidth]{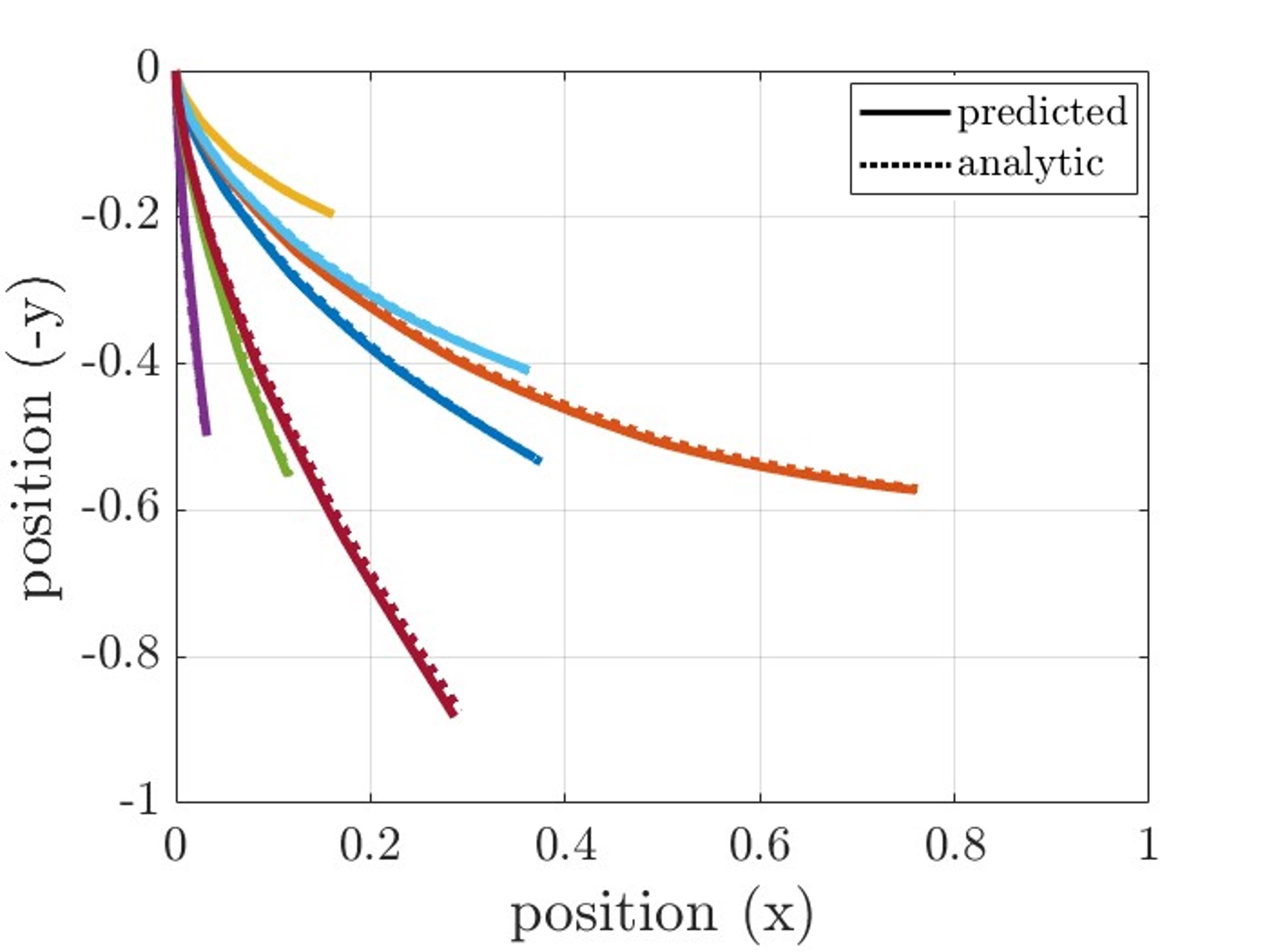}
    \caption{Solutions to the Brachistochrone problem.}
    \label{fig:brachi}
\end{figure}

Once trained, the Transformer was evaluated on a test dataset consisting of 10,000 trajectories to assess its performance and generalization capabilities. Table \ref{tab:brachi} shows that the Transformer predicts highly accurate results when compared to the analytical solutions, and can produce these results in real time. The table reports the following metrics: \textit{Traj. error} is the mean-squared-error between the results of the Transformer model and the analytical solution; \textit{Cost violation} is the average percent error between the evaluated cost of the predicted and analytical results; \textit{Loss} is the final loss value from Equation \eqref{eq:loss} for the Transformer; \textit{Computation time} is the average trajectory computation time from the Transformer model. Several trajectories are presented in Figure \ref{fig:brachi}, showing visually that the predicted output from the Transformer model is nearly identical to the analytical solutions.

%

\begin{table}
    \centering
    \caption{Prediction performance for Brachistochrone problem.} 
    \label{tab:brachi}
    \begin{tabular}{|l|c|c|}
        \hline
        \textbf{Metric} & \textbf{Transformer} \\ 
        \hline
        Traj. error (MSE)           & $6.13\times10^{-3}$           \\
        Cost violation (avg.)              & $1.35\%$              \\
        Loss (Eq. \ref{eq:loss}) & $1.13\times10^{-3}$  \\
        Computation time (avg.)        & $4.31\times10^{-2}\,\text{sec}$             \\
        \hline
    \end{tabular}
\end{table}

\begin{figure}[t]
    \centering
    \includegraphics[width=.48\textwidth]{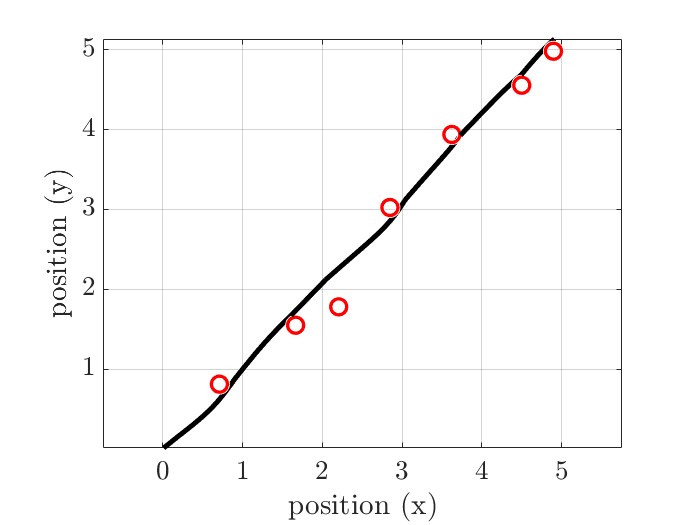}
    \caption{Sequential application of the obstacle avoidance Transformer model, avoiding multiple obstacles.}
    \label{fig:obsavoid}
\end{figure}

\begin{figure*}[t]
    \centering
    \begin{subfigure}{.32\textwidth}
        \centering
        \includegraphics[width=\textwidth]{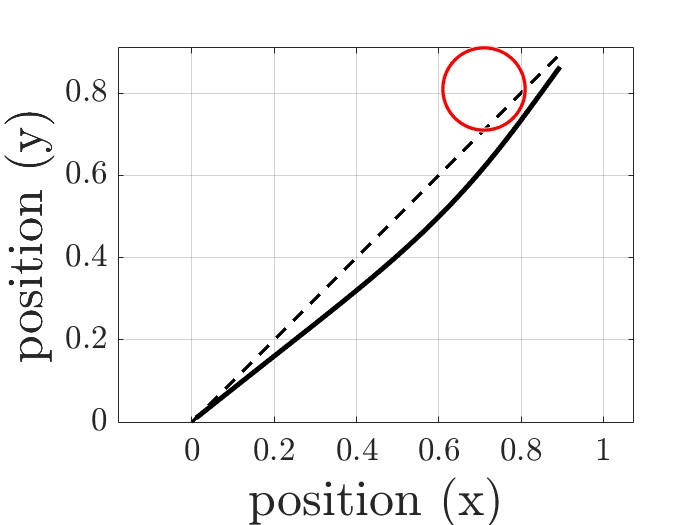}
        \caption{}\label{fig:obsavoid1}
    \end{subfigure}
    \begin{subfigure}{.32\textwidth}
        \centering
        \includegraphics[width=\textwidth]{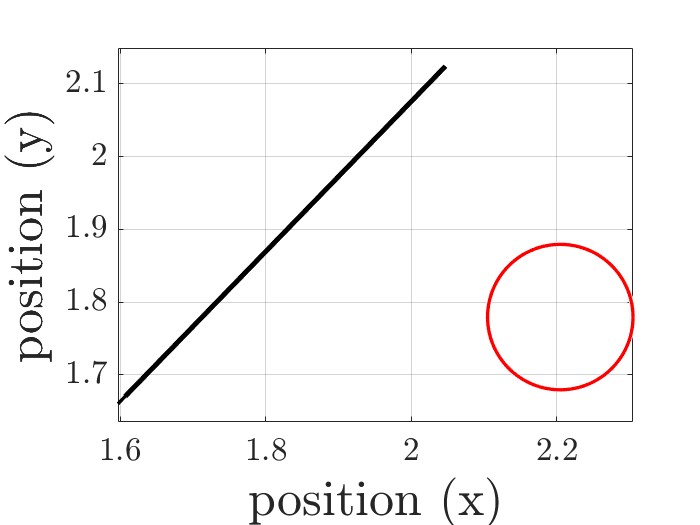}
        \caption{}\label{fig:obsavoid2}
    \end{subfigure}
    \begin{subfigure}{.32\textwidth}
        \centering
        \includegraphics[width=\textwidth]{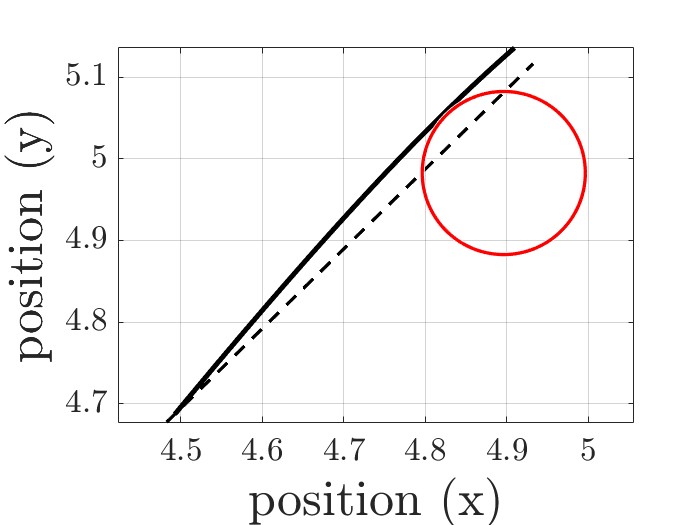}
        \caption{}\label{fig:obsavoid3}
    \end{subfigure}
    \caption{Three predictions from the obstacle avoidance Transformer model given some random obstacle position, with the dashed line being the trajectory in the absence of the obstacle.}
    \label{fig:obsavoid_sep}
\end{figure*}

\subsection{Obstacle avoidance problem}\label{subsec:obs}
We present an obstacle avoidance problem with single integrator dynamics:
\begin{equation}
    \min_{x_1(t),x_2(t),t_f}J(t_f)=t_f,
\end{equation}
subject to
\begin{equation}
\begin{split}
     \dot{x}_1(t)=u_1(t), \,\,\,\dot{x}_2(t)=u_2(t),
\end{split}
\end{equation}
\begin{equation}
\begin{split}
    & x_1(0)=\theta_1,  \,\,\,\,\,x_2(0)=\theta_2, \\
    & x_1(t_f)=\theta_3,\,\,\,x_2(t_f)=\theta_4,
\end{split}
\end{equation}
\begin{equation}
\begin{split}
    &u_{\text{min}}\le \sqrt{u_1^2(t)+u_2^2(t)}\le u_{\text{max}},\\
    &||x_1(t)-\theta_5||_2^2+||x_2(t)-\theta_6||^2_2\ge d,
\end{split}
\end{equation}
where $u_{\text{min}}$ and $u_{\text{max}}$ are velocity bounds, $\theta_1$, $\theta_2$, $\theta_3$, and $\theta_4$ are the initial and final $x$ and $y$ positions of the vehicle respectively, and $\theta_5$ and $\theta_6$ are the $x$ and $y$ positions of the obstacle. This problem can then be structured in the form of Problem $P^{OC}_M(\bm\theta)$ with $\bm{\theta}=[\theta_1,\theta_2,\theta_3,\theta_4,\theta_5,\theta_6]$, and then solved for $\bm{\theta}^i, \, i=1,...,$ 100,000.

This model can be used to generate potentially feasible, collision-free trajectories in real time in a receding-horizon fashion. At time $T$, let the position of the vehicle be $p_v(T)=[p_{v,x}(T), p_{v,y}(T)]$. Define $\theta_1=p_{v,x}(T)$ and $\theta_2=p_{v,y}(T)$. Let $\theta_3$ and $\theta_4$ represent a final position at a fixed distance from $p_v(T)$, oriented toward the desired destination. Finally, let $\theta_5$ and $\theta_6$ be the position of the obstacle seen by the vehicle, if any. Then, with $\bm{\theta}$ defined as above:
\begin{enumerate}
    \item Provide $\bm{\theta}$ to Algorithm \ref{alg:inference},
    \item Verify trajectory feasibility via degree elevation, convex hull, and end point properties; if infeasible, use as a warm-start to Problem $P_M^{OC}(\bm{\theta})$.
    \item Check if $\hat{\bar{\bm{x}}}_{M,M}$ is in the region of the desired terminal position, if so terminate, if not return to Step 1.
\end{enumerate}

Results from this planner are given in Figure \ref{fig:obsavoid}, which shows a vehicle avoiding a series of obstacles until reaching a desired terminal position. Three individual snapshot trajectories of Figure \ref{fig:obsavoid} are shown in Figure \ref{fig:obsavoid_sep}, where the solid line represents the predicted trajectory from the Transformer model and the dashed line indicates the optimal trajectory if there were no obstacle. As can be seen in Figures \ref{fig:obsavoid1} and \ref{fig:obsavoid3}, the Transformer is able to recognize that the nominal trajectory would result in collision, and successfully replans around the obstacle. Moreover, when the obstacle does not intersect with the nominal path, the Transformer model correctly predicts that the optimal solution lies on the straight line connecting the current and final position (see Figure \ref{fig:obsavoid2}).

\section{Conclusions}\label{sec:Conclusions}

In this work, we demonstrated the potential of Transformer models for prediction feasible solutions to infinite-dimensional optimization problems. By combining the ability of Transformers to capture long-range dependencies in sequential data with the geometric properties of composite Bernstein polynomials, we developed a framework of capable of generating feasible solutions online. Our results indicate that Transformer-based approaches can efficiently produce feasible trajectories for optimal control problems, highlighting their promise for real-time applications.

\bibliographystyle{IEEEtran}
\bibliography{refs}

@STRING{CDC  = {{IEEE Conference on Decision and Control}}}

@article{hauser2017learning,
  title={Learning the problem-optimum map: Analysis and application to global optimization in robotics},
  author={Hauser, Kris},
  journal={IEEE Transactions on Robotics},
  volume={33},
  number={1},
  pages={141--152},
  year={2017},
  publisher={IEEE}
}

@article{cichella2020optimal,
  title={Optimal multivehicle motion planning using {B}ernstein approximants},
  author={Cichella, Venanzio and Kaminer, Isaac and Walton, Claire and Hovakimyan, Naira and Pascoal, Antonio M},
  journal={IEEE transactions on automatic control},
  volume={66},
  number={4},
  pages={1453--1467},
  year={2020},
  publisher={IEEE}
}

@article{kielas2022bernstein,
  title={Bernstein Polynomial-Based Method for Solving Optimal Trajectory Generation Problems},
  author={Kielas-Jensen, Calvin and Cichella, Venanzio and Berry, Thomas and Kaminer, Isaac and Walton, Claire and Pascoal, Antonio},
  journal={Sensors},
  volume={22},
  number={5},
  pages={1869},
  year={2022},
  publisher={MDPI}
}

@inproceedings{cichella2019consistent,
  title={Consistent approximation of optimal control problems using {B}ernstein polynomials},
  author={Cichella, Venanzio and Kaminer, Isaac and Walton, Claire and Hovakimyan, Naira and Pascoal, Ant{\'o}nio M},
  booktitle={2019 IEEE 58th Conference on Decision and Control (CDC)},
  pages={4292--4297},
  year={2019},
  organization={IEEE}
}

@INPROCEEDINGS{maclin2024optimal,
  author={MacLin, Gage and Cichella, Venanzio and Patterson, Andrew and Acheson, Michael and Gregory, Irene},
  booktitle={2024 IEEE 63rd Conference on Decision and Control (CDC)}, 
  title={Optimal Control using Composite {B}ernstein Approximants}, 
  year={2024},
  volume={},
  number={},
  pages={8276-8281},
  doi={10.1109/CDC56724.2024.10886717}}

@article{vaswani2017attention,
  title={Attention is all you need},
  author={Vaswani, Ashish and Shazeer, Noam and Parmar, Niki and Uszkoreit, Jakob and Jones, Llion and Gomez, Aidan N and Kaiser, {\L}ukasz and Polosukhin, Illia},
  journal={Advances in neural information processing systems},
  volume={30},
  year={2017}
}

@article{cichella2017optimal,
  title={Optimal motion planning for differentially flat systems using {B}ernstein approximation},
  author={Cichella, Venanzio and Kaminer, Isaac and Walton, Claire and Hovakimyan, Naira},
  journal={IEEE Control Systems Letters},
  volume={2},
  year={2017},
  publisher={IEEE}
}

@article{tordesillas2022minvo,
  title={{MINVO} basis: Finding simplexes with minimum volume enclosing polynomial curves},
  author={Tordesillas, Jesus and How, Jonathan P},
  journal={Computer-Aided Design},
  volume={151},
  pages={103341},
  year={2022},
  publisher={Elsevier}
}

@inproceedings{guffanti2024transformers,
  title={Transformers for trajectory optimization with application to spacecraft rendezvous},
  author={Guffanti, Tommaso and Gammelli, Daniele and D’Amico, Simone and Pavone, Marco},
  booktitle={2024 IEEE Aerospace Conference},
  pages={1--13},
  year={2024},
  organization={IEEE}
}

@article{celestini2024transformer,
  title={Transformer-based model predictive control: Trajectory optimization via sequence modeling},
  author={Celestini, Davide and Gammelli, Daniele and Guffanti, Tommaso and D'Amico, Simone and Capello, Elisa and Pavone, Marco},
  journal={IEEE Robotics and Automation Letters},
  year={2024},
  publisher={IEEE}
}

@article{paszke2019pytorch,
  title={{P}y{T}orch: An Imperative Style, High-Performance Deep Learning Library},
  author={Paszke, Adam and Gross, Sam and Massa, Francisco and Lerer, Adam and Bradbury, James and Chanan, Gregory and Killeen, Trevor and Lin, Zeming and Gimelshein, Natalia and Antiga, Luca and others},
  journal={arXiv e-prints},
  pages={arXiv--1912},
  year={2019}
}

@article{malyuta2022convex,
  title={Convex optimization for trajectory generation: A tutorial on generating dynamically feasible trajectories reliably and efficiently},
  author={Malyuta, Danylo and Reynolds, Taylor P and Szmuk, Michael and Lew, Thomas and Bonalli, Riccardo and Pavone, Marco and A{\c{c}}{\i}kme{\c{s}}e, Beh{\c{c}}et},
  journal={IEEE Control Systems Magazine},
  volume={42},
  year={2022},
  publisher={IEEE}
}

@incollection{grune2016nonlinear,
  title={Nonlinear model predictive control: Theory and algorithms},
  author={Gr{\"u}ne, Lars and Pannek, J{\"u}rgen},
  pages={45--69},
  year={2016},
  publisher={Springer}
}

@book{liberzon2011calculus,
  title={Calculus of variations and optimal control theory: a concise introduction},
  author={Liberzon, Daniel},
  year={2011},
  publisher={Princeton university press}
}

@inproceedings{maclin2025transformer,
  title={Transformer-Enabled Leg-Based Trajectory Generation for Autonomous Aircraft Near Airports},
  author={MacLin, Gage and Cichella, Venanzio and Patterson, Andrew and Gregory, Irene M},
  booktitle={AIAA AVIATION FORUM AND ASCEND 2025},
  pages={3269},
  year={2025}
}

@article{patterson2023hermite,
  title={On {H}ermite Interpolation using {B}ernstein Polynomials for Trajectory Generation},
  author={Patterson, Andrew and MacLin, Gage and Acheson, Michael and Tabasso, Camilla and Cichella, Venanzio and Gregory, Irene},
  year={2023},
  booktitle={Tech. Rep. TM-2023-0013467, National Aeronautics and Space Administration}
}

\end{document}